\documentclass[12pt]{article}%
\usepackage{graphicx}
\usepackage[intlimits]{amsmath}
\usepackage{latexsym}
\usepackage{amsfonts}
\usepackage{amssymb}%
\makeatletter
\newfont{\footsc}{cmcsc10 at 8truept}
\newfont{\footbf}{cmbx10 at 8truept}
\newfont{\footrm}{cmr10 at 10truept}
\newtheorem{theorem}{Theorem}

\newtheorem{corollary}[theorem]{Corollary}

\newtheorem{problem}[theorem]{Problem}
\newtheorem{question}[theorem]{Question}

\newenvironment{proof}[1][Proof]{\noindent{\textbf {#1}  }}  {\hfill$\Box$\bigskip}

\begin{document}

\title{Graphs and Hermitian matrices: exact interlacing}
\author{B\'{e}la Bollob\'{a}s\thanks{Department of Mathematical Sciences, University
of Memphis, Memphis TN 38152, USA} \thanks{Department of Pure Mathematics \&
Mathematical Statistics University of Cambridge Centre for Mathematical
Sciences Wilberforce Road Cambridge CB3 0WB} \thanks{Research supported in
part by NSF grants CCR-0225610, DMS-0505550 and W911NF-06-1-0076.} \ and
Vladimir Nikiforov$^{\ast}$}
\maketitle

\begin{abstract}
We prove conditions for equality between the extreme eigenvalues of a matrix
and its quotient$.$ In particular, we give a lower bound on the largest
singular value of a matrix and generalize a result of Finck and Grohmann about
the largest eigenvalue of a graph.


\end{abstract}

\section{Introduction}

Our notation is standard (e.g., see \cite{Bol98}, \cite{CDS80}, and
\cite{HoJo88}); in particular, all graphs are defined on the vertex set
$\left[  n\right]  =\left\{  1,\ldots,n\right\}  $ and $G\left(  n\right)  $
stands for a graph of order $n$. Given a graph $G=G\left(  n\right)  ,$
$\mu_{1}\left(  G\right)  \geq...\geq\mu_{n}\left(  G\right)  $ are the
eigenvalues of its adjacency matrix $A\left(  G\right)  $, and $0=\lambda
_{1}\left(  G\right)  \leq...\leq\lambda_{n}\left(  G\right)  $ are the
eigenvalues of its Laplacian $L\left(  G\right)  $. If $X,Y\subset V\left(
G\right)  $ are disjoint sets, we write $G\left[  X\right]  $ for the graph
induced by $X,$ and $G\left[  X,Y\right]  $ for the bipartite graph induced by
$X$ and $Y;$ we set $e\left(  X\right)  =e\left(  G\left[  X\right]  \right)
$ and $e\left(  X,Y\right)  =e\left(  G\left[  X,Y\right]  \right)  $. We
assume that partitions consist of nonempty sets.

In this note we study conditions for finding exact eigenvalues using interlacing.

As proved in \cite{BoNi04}, if $G=G\left(  n\right)  $ and $\left[  n\right]
=\cup_{i=1}^{k}P_{i}$ is a partition, then%
\begin{equation}
\mu_{1}\left(  G\right)  +\ldots+\mu_{k}\left(  G\right)  \geq\sum_{i=1}%
^{k}\frac{2e\left(  P_{i}\right)  }{\left\vert P_{i}\right\vert },
\label{ineq4}%
\end{equation}%
\begin{equation}
\mu_{n-k+2}\left(  G\right)  +\ldots+\mu_{n}\left(  G\right)  \leq\sum
_{i=1}^{k}\frac{2e\left(  P_{i}\right)  }{\left\vert P_{i}\right\vert }%
-\frac{2e\left(  G\right)  }{n}, \label{ineq3}%
\end{equation}%
\begin{align}
\lambda_{2}\left(  G\right)  +\ldots+\lambda_{k}\left(  G\right)   &  \leq
\sum_{1\leq i<j\leq k}e\left(  P_{i},P_{j}\right)  \left(  \frac{1}{\left\vert
P_{i}\right\vert }+\frac{1}{\left\vert P_{j}\right\vert }\right)
,\label{lapl1}\\
\lambda_{n-k+1}\left(  G\right)  +\ldots+\lambda_{n}\left(  G\right)   &
\geq\sum_{1\leq i<j\leq k}e\left(  P_{i},P_{j}\right)  \left(  \frac
{1}{\left\vert P_{i}\right\vert }+\frac{1}{\left\vert P_{j}\right\vert
}\right)  . \label{lapl2}%
\end{align}

To warm up we shall give necessary conditions for equality in these inequalities.

Call a bipartite graph \emph{semiregular} if the vertices of the same vertex
class have the same degree. Call a partition $V\left(  G\right)  =\cup
_{i=1}^{k}P_{i}$ \emph{semiequitable for }$G$ if $G\left[  P_{i},P_{j}\right]
$ is semiregular for $1\leq i<j\leq k,$ and \emph{equitable for }$G$ if, in
addition, $G\left[  P_{i}\right]  $ is regular for $i\in\left[  k\right]  $.

\begin{theorem}
\label{th1} If equality holds in (\ref{ineq4}) or (\ref{ineq3}), then the
partition $\left[  n\right]  =\cup_{i=1}^{k}P_{i}$ is equitable for $G$;
moreover, if equality holds in (\ref{ineq3}), then $G$ is regular. If equality
holds in (\ref{lapl1}) or (\ref{lapl2}), then the partition $\left[  n\right]
=\cup_{i=1}^{k}P_{i}$ is semiequitable for $G.$
\end{theorem}

In order to discuss this result from a more general viewpoint, we introduce
additional notation and definitions. We order the eigenvalues of an $n\times
n$ Hermitian matrix $A$ as $\mu_{1}\left(  A\right)  \geq...\geq\mu_{n}\left(
A\right)  .$

Suppose $1<k<n$ and let $A$ and $B$ be Hermitian matrices of size $n\times n$
and $k\times k$. As usual, we say that the eigenvalues of $A$ and $B$ are
\emph{interlaced}, if $\mu_{i}\left(  A\right)  \geq\mu_{i}\left(  B\right)
\geq\mu_{n-k+i}\left(  A\right)  $ for all $i\in\left[  k\right]  $. The
interlacing is called \emph{tight} if there exists an integer $r\in\left[
0,k\right]  $ such that%
\[
\mu_{i}\left(  A\right)  =\mu_{i}\left(  B\right)  \text{ for }0\leq i\leq
r\text{ and }\mu_{n-k+i}\left(  A\right)  =\mu_{i}\left(  B\right)  \text{ for
}r<i\leq k.
\]
When we must indicate the value $r,$ we say that the interlacing is $r$-tight.

Note that inequalities (\ref{ineq4}) - (\ref{lapl2}) are proved using
eigenvalue interlacing; we shall see that equality in either of them implies
tight interlacing, in turn, implying the conditions of Theorem \ref{th1}.
Hence, the following question arises.

\begin{question}
\textit{\label{pro1}For which graphs conditions similar to those in Theorem
\ref{th1} imply tight interlacing.}
\end{question}

Below we answer a simple, yet important case of this question.

\begin{theorem}
\label{th2}Let $G=G\left(  n\right)  $ and $\left[  n\right]  =\cup_{i=1}%
^{k}P_{i}$ be a partition such that, for all $i\in\left[  k\right]  $,
$G\left[  P_{i}\right]  $ is empty and, for all $1\leq i<j\leq k$, $G\left[
P_{i},P_{j}\right]  $ is empty or complete. Then equality holds in
(\ref{ineq4}), (\ref{lapl1}), and (\ref{lapl2}). If $G$ is regular, equality
holds in (\ref{ineq3}) as well.
\end{theorem}

The general case of Question \ref{pro1} seems rather difficult; however, most
often we are interested in simpler problems, which, for convenience, we state
for matrices.

We first relax the concept of tight interlacing. Suppose $1<k<n$ and let $A$
and $B$ be Hermitian matrices of size $n\times n$ and $k\times k$ with
interlaced eigenvalues. Call the interlacing \emph{exact} if there exist
integers $p,q$ such that $0<p+q\leq k$ and%
\[
\mu_{i}\left(  A\right)  =\mu_{i}\left(  B\right)  \text{ for }0\leq i\leq
p\text{ and }\mu_{n-k+i}\left(  A\right)  =\mu_{i}\left(  B\right)  \text{ for
}k-q<i\leq k.
\]
When we must indicate the values $p$ and $q,$ we say that the interlacing is
$\left(  p,q\right)  $-exact.

\begin{problem}
\label{pro2} Find conditions for $\left(  p,q\right)  $-exact interlacing.
\end{problem}

Among all combinations of $p$ and $q,$ the case of $p=1,$ $q=0$ is of primary
importance. We give a solution to Problem \ref{pro2} in this case, when $A$
is\ nonnegative and $B$ is a \textquotedblleft quotient\textquotedblright%
\ matrix of $A$. Again, we introduce some notation and definitions.

Given an $m\times n$ matrix $A=\left\{  a_{ij}\right\}  $ and sets
$I\subset\left[  m\right]  ,$ $J\subset\left[  n\right]  ,$ write $A\left[
I,J\right]  $ for the submatrix of all $a_{ij}$ with $i\in I$ and $j\in J.$ A
matrix $A$ is called \emph{regular} if its row sums are equal and so are its
column sums.

Let $A=\left\{  a_{ij}\right\}  $ be an $m\times n$ matrix and let
$\mathcal{P}=\left\{  P_{1},\ldots,P_{k}\right\}  $, $\mathcal{Q}=\left\{
Q_{1},\ldots,Q_{l}\right\}  $ be partitions of $\left[  m\right]  $ and
$\left[  n\right]  .$ Set $\mathcal{P}\times\mathcal{Q}=\left\{  P_{i}\times
Q_{j}:i\in\left[  k\right]  ,j\in\left[  l\right]  \right\}  $ and note that
$\mathcal{P}\times\mathcal{Q}$ is a partition of $\left[  m\right]
\times\left[  n\right]  $. Call the partition $\mathcal{P}\times\mathcal{Q}$
\emph{equitable for }$A$ if $A\left[  P_{p},Q_{q}\right]  $ is regular for all
$p\in\left[  k\right]  ,$ $q\in\left[  l\right]  .$ Write $A|\mathcal{P}%
\times\mathcal{Q}$ for the $k\times l$ matrix $\left\{  b_{pq}\right\}  $
defined by
\[
b_{pq}=\frac{1}{\sqrt{\left\vert P_{p}\right\vert \left\vert Q_{q}\right\vert
}}\sum_{i\in P_{p},j\in Q_{q}}a_{ij}\text{, \ \ \ }p\in\left[  k\right]
,q\in\left[  l\right]  .
\]
Sometimes $A|\mathcal{P}\times\mathcal{Q}$ is called a \emph{quotient matrix}
of $A.$

Haemers \cite{Hae95} proved the following result.

\begin{theorem}
\label{ThH}For any $n\times n$ Hermitian matrix $A$ and any partition
$\mathcal{P}$ of $\left[  n\right]  ,$ the eigenvalues of $A$ and
$A|\mathcal{P}\times\mathcal{P}$ are interlaced; moreover, if the interlacing
is tight then $\mathcal{P}\times\mathcal{P}$ is equitable for $A$.
\end{theorem}

In particular, for any $n\times n$ Hermitian matrix $A$ and any partition
$\mathcal{P}$ of $\left[  n\right]  ,$ $\mu_{1}\left(  A\right)  \geq\mu
_{1}\left(  A|\mathcal{P}\times\mathcal{P}\right)  .$ We use the
Perron-Frobenius theorem to prove sufficient conditions for equality in this inequality.

\begin{theorem}
\label{th3}If $A$ is an irreducible, nonnegative symmetric matrix and
$\mathcal{P}\mathbf{\times}\mathcal{P}$ is equitable for $A$, then $\mu
_{1}\left(  A\right)  =\mu_{1}\left(  A|\mathcal{P}\times\mathcal{P}\right)  $.
\end{theorem}

We deduce a similar result about the largest singular value of a matrix. Write
$A^{\ast}$ for the Hermitian transpose of $A.$

\begin{theorem}
\label{th4}Let $A$ be a complex $m\times n$ matrix, $\mathcal{P}$ a partition
of $\left[  m\right]  $, and $\mathcal{Q}$ a partition of $\left[  n\right]
.$ Then $\sigma_{1}\left(  A\right)  \geq\sigma_{1}\left(  A|\mathcal{P}%
\times\mathcal{Q}\right)  .$

If $A$ is nonnegative, $AA^{\ast}$ and $A^{\ast}A$ are irreducible, and
$\mathcal{P}\mathbf{\times}\mathcal{Q}$ is equitable for $A$, then $\sigma
_{1}\left(  A\right)  =\sigma_{1}\left(  A|\mathcal{P}\times\mathcal{Q}%
\right)  $.
\end{theorem}

Note that the first part of this result is implicit in \cite{Hae95}. Observe
also that the conditions for equality in Theorem \ref{th3} and Theorem
\ref{th4} are sufficient but not necessary. Thus, we have another question.

\begin{question}
\label{pro3}For which nonnegative $m\times n$ matrices $A$ and partitions
$\mathcal{P}$ of $\left[  m\right]  $, and $\mathcal{Q}$ of $\left[  n\right]
,$ does the condition that $\mathcal{P}\times\mathcal{Q}$ is equitable for
$A$\ imply that $\sigma_{1}\left(  A\right)  =\sigma_{1}\left(  A|\mathcal{P}%
\times\mathcal{Q}\right)  $?
\end{question}

We can answer Question \ref{pro3} in a particular case, generalizing a
classical result on graph spectra. Write $G_{1}+G_{2}$ for the join of the
graphs $G_{1}$, $G_{2}$ and recall a theorem of Finck and Grohmann
\cite{FiGr65} (see also \cite{CDS80}, Theorem 2.8):

\textit{Let the graph }$G$ be the join of an $r_{1}$\textit{-regular graph
}$G_{1}$ of order $n_{1}$\textit{ and an }$r_{2}$-regular \textit{graph
}$G_{2}$ of order $n_{2}.$\textit{ Then }$\mu_{1}\left(  G\right)  $\textit{
is the positive root of the equation }%
\begin{equation}
\left(  x-r_{1}\right)  \left(  x-r_{2}\right)  -n_{1}n_{2}=0. \label{eq1}%
\end{equation}

Setting $\mathcal{P}=\left\{  V\left(  G_{1}\right)  ,V\left(  G_{2}\right)
\right\}  ,$ a routine calculation shows that (\ref{eq1}) is the
characteristic equation of $A\left(  G\right)  |\mathcal{P}\times\mathcal{P}$;
therefore, the conclusion of the Finck-Grohmann theorem reads as
\[
\mu_{1}\left(  G\right)  =\mu_{1}\left(  A\left(  G\right)  |\mathcal{P}%
\times\mathcal{P}\right)  .
\]

Clearly if $G_{i}=G\left(  n_{i}\right)  ,$ $2\leq i\leq k,$ and
$G=G_{1}+\ldots+G_{k},$ then letting $\mathcal{P}=\left\{  V\left(
G_{1}\right)  ,\ldots,V\left(  G_{k}\right)  \right\}  ,$ by Theorem
\ref{ThH},
\[
\mu_{1}\left(  G\right)  \geq\mu_{1}\left(  A\left(  G\right)  |\mathcal{P}%
\times\mathcal{P}\right)  .
\]

It is natural to ask when $\mu_{1}\left(  G\right)  =\mu_{1}\left(  A\left(
G\right)  |\mathcal{P}\times\mathcal{P}\right)  .$ We deduce the answer of
this question from a more general matrix result.

\begin{theorem}
\label{th5}Let $A$ be a symmetric, irreducible, nonnegative matrix of size
$n\times n$ and $\mathcal{P}=\left\{  P_{1},\ldots,P_{k}\right\}  $ be a
partition of $\left[  n\right]  $ such that $A\left[  P_{i},P_{j}\right]  $ is
regular for all $1\leq i<j\leq k$. Then
\begin{equation}
\mu_{1}\left(  A\right)  =\mu_{1}\left(  A|\mathcal{P}\times\mathcal{P}%
\right)  \label{in5}%
\end{equation}
if and only if $A\left[  P_{i},P_{i}\right]  $ is regular for all $1\leq i\leq
k,$ i.e., $\mathcal{P}\times\mathcal{P}$ is regular in $A.$
\end{theorem}

For graphs this theorem implies the following corollary.

\begin{corollary}
\label{co1}Let $G=G\left(  n\right)  $ be a connected graph and $\mathcal{P}$
be a semiequitable for $G$ partition of $\left[  n\right]  $. Then $\mu
_{1}\left(  G\right)  =\mu_{1}\left(  A\left(  G\right)  |\mathcal{P}%
\times\mathcal{P}\right)  $ if and only if $\mathcal{P}$ is equitable for $G.$
\end{corollary}

\section{Proofs}

\begin{proof}
[\textbf{Proof of Theorem \ref{th1}}]For short set $A=A\left(  G\right)  $ and
$L=L\left(  G\right)  .$ Equality in (\ref{ineq4}) implies that
\[
\mu_{1}\left(  G\right)  +\ldots+\mu_{k}\left(  G\right)  =\sum_{i=1}^{k}%
\frac{2e\left(  P_{i}\right)  }{\left\vert P_{i}\right\vert }=tr\left(
A|\mathcal{P}\times\mathcal{P}\right)  ;
\]
hence $\mu_{i}\left(  G\right)  =\mu_{i}\left(  A|\mathcal{P}\times
\mathcal{P}\right)  $ for all $i\in\left[  k\right]  .$ Thus, the interlacing
is $k$-tight and $\mathcal{P}\times\mathcal{P}$ is equitable for $A$:
therefore, $\mathcal{P}$ is equitable for $G.$

Inequality (\ref{ineq3}) follows from Theorem \ref{ThH} and $\mu_{1}\left(
A|\mathcal{P}\times\mathcal{P}\right)  \geq2e\left(  G\right)  /n$ noting
that
\begin{align*}
\mu_{n-k+2}\left(  G\right)  +\ldots+\mu_{n}\left(  G\right)   &  \leq
tr\left(  A|\mathcal{P}\times\mathcal{P}\right)  -\mu_{1}\left(
A|\mathcal{P}\times\mathcal{P}\right) \\
&  =\sum_{i=1}^{k}\frac{2e\left(  P_{i}\right)  }{\left\vert P_{i}\right\vert
}-\mu_{1}\left(  A|\mathcal{P}\times\mathcal{P}\right)  .
\end{align*}
Hence, if equality holds in (\ref{ineq3}), then $\mu_{n-k+i}\left(  G\right)
=\mu_{i}\left(  A|\mathcal{P}\times\mathcal{P}\right)  $ for every
$i=2,\ldots,k.$ To prove that the interlacing is tight, we shall show that
$\mu_{1}\left(  G\right)  =\mu_{1}\left(  A|\mathcal{P}\times\mathcal{P}%
\right)  .$ Note first $\mu_{1}\left(  A|\mathcal{P}\times\mathcal{P}\right)
=2e\left(  G\right)  /n$. Also it is easy to see that the $k$-vector $\left(
\sqrt{\left\vert P_{1}\right\vert },\ldots,\sqrt{\left\vert P_{k}\right\vert
}\right)  $ is an eigenvector to $\mu_{1}\left(  A|\mathcal{P}\times
\mathcal{P}\right)  .$ This implies that the $n$-vector of all ones is an
eigenvector to $G$ and $\mu_{1}\left(  A|\mathcal{P}\times\mathcal{P}\right)
$ is an eigenvalue of $G;$ hence, the Perron-Frobenius theorem implies that
$G$ is regular and $\mu_{1}\left(  G\right)  =2e\left(  G\right)  /n=\mu
_{1}\left(  A|\mathcal{P}\times\mathcal{P}\right)  .$ Therefore, the
interlacing is $1$-tight and $\mathcal{P}\times\mathcal{P}$ is equitable for
$A;$ so $\mathcal{P}$ is equitable for $G.$

Inequality (\ref{lapl1}) follows from Theorem \ref{ThH} and $\mu_{k}\left(
L|\mathcal{P}\times\mathcal{P}\right)  =0$ noting that
\begin{align*}
\lambda_{2}\left(  G\right)  +\ldots+\lambda_{k}\left(  G\right)   &  \leq
\sum_{i=1}^{k-1}\mu_{i}\left(  L|\mathcal{P}\times\mathcal{P}\right)
=tr\left(  L|\mathcal{P}\times\mathcal{P}\right) \\
&  \mathbf{=}\sum_{1\leq i<j\leq k}e\left(  P_{i},P_{j}\right)  \left(
\frac{1}{\left\vert P_{i}\right\vert }+\frac{1}{\left\vert P_{j}\right\vert
}\right)  .
\end{align*}
Consequently, by $\lambda_{1}\left(  G\right)  =0,$ equality in (\ref{lapl1})
implies that the interlacing is $1$-tight. Hence, $\mathcal{P}\times
\mathcal{P}$ is equitable for $L$, and so, for all $1\leq i<j\leq k,$ the
graphs $G\left[  P_{i},P_{j}\right]  $ are semiregular.

Finally, inequality (\ref{lapl2}) follows from Theorem \ref{ThH} noting that%
\[
\lambda_{n-k+1}\left(  G\right)  +\ldots+\lambda_{n}\left(  G\right)  \geq
\sum_{i=1}^{k}\mu_{i}\left(  L|\mathcal{P}\times\mathcal{P}\right)
=\sum_{1\leq i<j\leq k}e\left(  P_{i},P_{j}\right)  \left(  \frac
{1}{\left\vert P_{i}\right\vert }+\frac{1}{\left\vert P_{j}\right\vert
}\right)  .
\]
Clearly, equality in (\ref{lapl2}) implies that the interlacing is $\left(
k-1\right)  $-tight. Hence, $\mathcal{P}\times\mathcal{P}$ is equitable for
$L;$ thus, for all $1\leq i<j\leq k,$ the graphs $G\left[  P_{i},P_{j}\right]
$ are semiregular, as claimed.
\end{proof}

\bigskip

\begin{proof}
[\textbf{Proof of Theorem \ref{th2}}]For short write $A$ for $A\left(
G\right)  .$ Since $\mathcal{P}\times\mathcal{P}$ is equitable for $A,$ for
every unit eigenvector $\mathbf{y}=\left(  y_{1},\ldots,y_{k}\right)  $ to an
eigenvalue $\mu$ of $A|\mathcal{P}\times\mathcal{P}$, the vector
$\mathbf{x}=\left(  x_{1},\ldots,x_{n}\right)  $ defined by%
\[
x_{i}=\frac{1}{\sqrt{\left\vert P_{s}\right\vert }}y_{s}\text{ \ \ for \ }i\in
P_{s}%
\]
is a unit eigenvector of $A$ to the eigenvalue $\mu.$ This implies that the
spectrum of $A$ contains all eigenvalues of $A|\mathcal{P}\times\mathcal{P}$
with the same or greater multiplicity.

On the other hand, the structure of $G$ implies that the vertices in the same
partition set $P_{i}$ have the same neighbors. Thus, every eigenvalue $\mu$ of
$A$ has an eigenvector which is constant within each $P_{i}$. This implies
that every eigenvalue of $A$ is also an eigenvalue of $A|\mathcal{P}%
\times\mathcal{P}$. Therefore, $A$ and $A|\mathcal{P}\times\mathcal{P}$ have
the same set of eigenvalues and each eigenvalue occurs at least as many times
in the spectrum of $A$ as in the spectrum of $A|\mathcal{P}\times\mathcal{P}$.
Since $tr\left(  A^{2}\right)  =tr\left(  A|\mathcal{P}\times\mathcal{P}%
\right)  ^{2},$ we see that
\[
\sum_{i=1}^{n}\mu_{i}^{2}\left(  A\right)  =\sum_{i=1}^{n}\mu_{i}^{2}\left(
A|\mathcal{P}\times\mathcal{P}\right)  ,
\]
and so $A|\mathcal{P}\times\mathcal{P}$ and $A$ have exactly the same nonzero
eigenvalues with the same multiplicities. Hence, inequalities (\ref{ineq4}%
)-(\ref{lapl2}) follow immediately, completing the proof.
\end{proof}

\bigskip

\begin{proof}
[\textbf{Proof of Theorem \ref{th3}}]Let $\mathcal{P}=\left\{  P_{1}%
,\ldots,P_{k}\right\}  $ and suppose that $\mathcal{P}\times\mathcal{P}$ is
equitable for $A.$ Since $A$ is irreducible, $A|\mathcal{P}\times\mathcal{P}$
is also irreducible; let $\mathbf{y}=\left(  y_{1},\ldots,y_{k}\right)  $ be a
positive unit eigenvector to $\mu_{1}\left(  A|\mathcal{P}\times
\mathcal{P}\right)  .$ Then the vector $\mathbf{x}=\left(  x_{1},\ldots
,x_{n}\right)  $ defined by%
\[
x_{i}=\frac{1}{\sqrt{\left\vert P_{s}\right\vert }}y_{s}\text{ \ \ for
\ \ }i\in P_{s}%
\]
is a positive unit vector such that $A\mathbf{x}=\mu_{1}\left(  A|\mathcal{P}%
\times\mathcal{P}\right)  \mathbf{x},$ implying that $\mu_{1}\left(
A|\mathcal{P}\times\mathcal{P}\right)  $ is an eigenvalue of $A$ with
eigenvector $\mathbf{x}$. The Perron-Frobenius theorem implies that $\mu
_{1}\left(  A|\mathcal{P}\times\mathcal{P}\right)  =\mu_{1}\left(  A\right)
,$ completing the proof.
\end{proof}

\bigskip

\begin{proof}
[\textbf{Proof of Theorem \ref{th4}}]For every $i\in\left[  l\right]  $, set
$Q_{i}^{\prime}=\left\{  x+m:x\in Q_{i}\right\}  $; thus $\mathcal{Q}^{\prime
}=\left\{  Q_{1}^{\prime},\ldots,Q_{l}^{\prime}\right\}  $ is a partition of
the set $\left[  m+1,m+n\right]  $ and $\mathcal{R}\mathbf{=}\mathcal{P}%
\cup\mathcal{Q}^{\prime}$ is a partition of $\left[  m+n\right]  .$ Let
\[
B=\left(
\begin{array}
[c]{cc}%
0 & A^{\ast}\\
A & 0
\end{array}
\right)  .
\]
It is known (see, e.g., \cite{HoJo88}, p. 418) that $\sigma_{1}\left(
A\right)  =\mu_{1}\left(  B\right)  .$ It is easy to see that $B$ is
irreducible if and only if $A^{\ast}A$ and $AA^{\ast}$ are irreducible. Since
\[
B|\mathcal{R}\times\mathcal{R}=\left(
\begin{array}
[c]{cc}%
0 & \left(  A|\mathcal{P}\times\mathcal{Q}\right)  ^{\ast}\\
A|\mathcal{P}\times\mathcal{Q} & 0
\end{array}
\right)  ,
\]
if $\sigma_{1}\left(  A\right)  =\sigma_{1}\left(  A|\mathcal{P}%
\times\mathcal{Q}\right)  ,$ we see that
\[
\mu_{1}\left(  B\right)  =\sigma_{1}\left(  A\right)  =\sigma_{1}\left(
A|\mathcal{P}\times\mathcal{Q}\right)  =\mu_{1}\left(  B|\mathbf{R}%
\times\mathbf{R}\right)  ,
\]
and Theorem \ref{th3} implies that $\mathcal{R}\times\mathcal{R}$ is equitable
for $B;$ hence, $\mathcal{P}\times\mathcal{Q}$ is equitable for $A$,
completing the proof.
\end{proof}

\bigskip

\begin{proof}
[\textbf{Proof of Theorem \ref{th5}}]If $A\left[  P_{i},P_{i}\right]  $ is
regular for each $i\in\left[  k\right]  ,$ the partition $\mathcal{P}%
\times\mathcal{P}$ is equitable for $A,$ and Theorem \ref{th3} implies
(\ref{in5}). Suppose now $\mu_{1}\left(  A\right)  =\mu_{1}\left(
A|\mathcal{P}\times\mathcal{P}\right)  .$ We have to prove that $A\left[
P_{i},P_{i}\right]  $ is regular for every $i\in\left[  k\right]  .$ Since
$A|\mathcal{P}\times\mathcal{P}$ is irreducible, there is a positive unit
eigenvector $\mathbf{y}=\left(  y_{1},\ldots,y_{k}\right)  $ to $\mu
_{1}\left(  A|\mathcal{P}\times\mathcal{P}\right)  $. Define the unit vector
$\mathbf{x}=\left(  x_{1},\ldots,x_{n}\right)  $ by%
\[
x_{i}=\frac{1}{\sqrt{\left\vert P_{s}\right\vert }}y_{s}\text{ \ \ for
\ \ }i\in P_{s}.
\]
We have $\left\langle A\mathbf{x},\mathbf{x}\right\rangle =\mu_{1}\left(
A|\mathcal{P}\times\mathcal{P}\right)  =\mu_{1}\left(  A\right)  ,$ and so
$\mathbf{x}$ is an eigenvector of $A$ to $\mu_{1}\left(  A\right)  .$ For any
$r\in\left[  k\right]  $ and $s,t\in P_{r},$ we have
\begin{align*}
\mu_{1}\left(  A\right)  x_{s}  &  =\sum_{i=1}^{n}a_{si}x_{i}=\sum_{i=1}%
^{k}\frac{1}{\sqrt{\left\vert P_{i}\right\vert }}y_{i}\sum_{j\in P_{i}}%
a_{sj}\\
\mu_{1}\left(  A\right)  x_{t}  &  =\sum_{i=1}^{n}a_{ti}x_{i}=\sum_{i=1}%
^{k}\frac{1}{\sqrt{\left\vert P_{i}\right\vert }}y_{i}\sum_{j\in P_{i}}a_{tj}%
\end{align*}
Since $x_{s}=x_{t}$ and
\[
\sum_{j\in P_{i}}a_{tj}=\sum_{j\in P_{i}}a_{sj}%
\]
for $i\in\left[  k\right]  \backslash\left\{  r\right\}  ,$ we see that
\[
\sum_{j\in P_{r}}a_{sj}=\sum_{j\in P_{r}}a_{tj},
\]
that is to say, the row sums of $A\left[  P_{r},P_{r}\right]  $ are equal.
Since $A\left[  P_{r},P_{r}\right]  $ is symmetric, this implies that
$A\left[  P_{r},P_{r}\right]  $ is regular, completing the proof.
\end{proof}

\bigskip

\textbf{Concluding remarks}

In this note we confined our investigation of exact interlacing to the largest
eigenvalue only. It would be good to continue this work for the smallest and
the second largest eigenvalues, i.e., for $\left(  0,1\right)  $-exact and
$\left(  2,0\right)  $-exact interlacing. Unfortunately, these important
problems seem rather difficult to tackle.

\end{document}